\documentclass[12pt,a4paper]{article}

\usepackage{pdflscape}
\usepackage{amsmath}
\usepackage{amssymb}
\usepackage{latexsym}
\usepackage{srcltx}
\usepackage{graphics}
\usepackage{color}
\usepackage{epsfig}
\usepackage{color}
\usepackage{url}

\newcommand{\co}{{\mathbb C}}

\newcommand{\re}{{\mathbb R}}
\newcommand{\n}{{\mathbb N}}

\newcommand{\cA}{{\mathcal{A}}}

\newcommand{\cV}{{\mathcal{V}}}

\newcommand{\cB}{{\mathcal{B}}}

\newcommand{\cS}{{\mathcal{S}}}
\newcommand{\cR}{{\mathcal{R}}}
\newcommand{\cQ}{{\mathcal{Q}}}

\newcommand{\bE}{{\boldsymbol{E}}}
\newcommand{\bH}{{\boldsymbol{H}}}
\newcommand{\bK}{{\boldsymbol{K}}}
\newcommand{\bL}{{\boldsymbol{L}}}
\newcommand{\bM}{{\boldsymbol{M}}}
\newcommand{\bR}{{\boldsymbol{R}}}

\newcommand{\bU}{{\boldsymbol{U}}}
\newcommand{\bV}{{\boldsymbol{V}}}

\newcommand{\bx}{{\boldsymbol{x}}}
\newcommand{\by}{{\boldsymbol{y}}}

\newcommand{\be}{{\boldsymbol{e}}}

\newcommand{\ba}{{\boldsymbol{a}}}
\newcommand{\bb}{{\boldsymbol{b}}}

\newcommand{\bv}{{\boldsymbol{v}}}
\newcommand{\bu}{{\boldsymbol{u}}}

\newcommand{\bg}{{\boldsymbol{g}}}

\newcommand{\nill}{{\boldsymbol{0}}}

\newtheorem{theorem}{Theorem}
\newtheorem{prop}{Proposition}
\newtheorem{lemma}{Lemma}
\newtheorem{cor}{Corollary}
\newtheorem{remark}{Remark}

\newtheorem{defi}{Definition}

\newtheorem{assum}{Assumption}

\usepackage{soul}
\usepackage{cancel}
\usepackage[normalem]{ulem}

\usepackage[dvipsnames]{xcolor}
\usepackage{tikz}
\usepackage{tkz-euclide}

\usepackage{rotating}

\usepackage{multicol}
\usepackage[hidelinks]{hyperref}

\date{}

\medskip 
\author{
Vladimir Yu. Protasov 
\thanks{University of L'Aquila, Italy; {e-mail: \tt\small
vladimir.protasov@univaq.it}} 
}

\title{Perron matrix semigroups}

\begin{document}
\maketitle

\begin{abstract}

We consider multiplicative semigroups of real~$d\times d$ matrices. A semigroup~$\cS$ is 
called Perron if each of its matrices has a Perron eigenvalue, i.e., an eigenvalue 
equal to the spectral radius. If all matrices of~$\cS$ 
leave a  proper convex cone invariant, then~$\cS$ is Perron. 
Our main result asserts   the converse: every irreducible Perron semigroup 
possesses a common invariant cone,   provided that some mild 
assumptions are satisfied.  This gives conditions for a set of matrices to 
share a common invariant cone, which is an important property widely studied in the literature. 
Then we address the problem to characterise the exceptions, when a Perron semigroup 
does not have an invariant cone.  For~$d\le 4$, all Perron semigroups are classified.  
For higher dimensions~$d$,  several classes of such semigroups are found.

\bigskip

\smallskip

\begin{flushright}

\end{flushright}

\end{abstract}
\bigskip

\begin{center}
\large{\textbf{1. Introduction}}
\end{center}
\bigskip 

\vspace{1cm}

If several real~$d\times d$ matrices have a common invariant cone, 
then they possess  special properties  similar to properties of non-negative
matrices. In particular, some of the key facts of the Perron-Frobenius theory  can be generalized to 
matrix families  sharing a common  invariant cone. A question arises how  
 to prove those properties if an invariant cone is not given? 
 And how to determine that such a cone exists at all?  The problem to  decide whether a given set of matrices 
 possesses  a  common invariant cone is referred to as 
the {\em invariant cone problem}. It  is notoriously hard. A short overview  is given in Section~2.  In this paper, we obtain necessary and sufficient conditions.  They are  close to each other and can thus be considered 
as a criterion  which holds apart from some exceptional cases. 
Classification of those exceptions is an  intriguing problem related to matrix semigroups with 
special spectral properties. This issue is 
addressed in Sections 6 and 7.
\smallskip

Consider an arbitrary nontrivial (containing at least one nonzero element) multiplicative  semigroup of real  $d\times d$  matrices~$\cS$. 
This means that for every~$A, B\in \cS$, we have~$AB\in \cS$. 
Let~$\bK\subset \re^d$ be a closed convex cone with an apex at the origin, i.e., a closed  set such that for every 
$\bx, \by \in \bK$, we have $\bx + \by \in \bK$ and for every number~$t \ge 0$, we have~$t\bx \in \bK$. The cone is called {\em proper} if it is  non-degenerate (possesses a nonempty interior) and pointed, i.e., $\bK \cap (-\bK) \, = \, \{0\}$. 
All cones are assumed to be proper. 
We say that a  matrix~$A$ leaves a cone $\bK$ invariant if~$A\bK \subset \bK$. 
If all matrices from a semigroup~$\cS$ leave~$\bK$ invariant, 
then~$\bK$ is called an {\em invariant cone} of~$\cS$.   

For an arbitrary $d\times d$ matrix~$A$, its maximal by modulus eigenvalues are referred to as {\em leading}. If a leading eigenvalue is nonnegative, it is called~{\em Perron}. 
By the Krein-Rutman theorem~\cite{KR}, a matrix that leaves a cone
invariant  possesses a Perron eigenvalue. This implies that not every matrix 
has an invariant cone.  Moreover, the existence of Perron eigenvalue is, in general, not sufficient for that.  A criterion for a matrix to possess an invariant cone is presented in~\cite{V}. 
\begin{defi}\label{d.5}
A multiplicative matrix semigroup~$\cS$ is called a {\em Perron semigroup} 
if every element of~$\cS$ possesses a Perron eigenvalue. 
\end{defi}
An immediate conclusion from the Krein-Rutman theorem is the following: 
\medskip 
 
\noindent \textbf{Fact~1.} {\em Every  semigroup that has an invariant cone is Perron}. 
\medskip 

  Our interest will be in establishing  the 
 converse. This could give a criterion for an arbitrary set of matrices to possess an  invariant cone. Under which conditions, every Perron semigroup necessarily has an invariant cone? 
 Let us stress that we require the existence of a {\em common invariant cone}
 for all matrices. For example, if all matrices from some finite set~$\cA = \{A_1, \ldots , A_m\}$ and all their products possess 
 Perron eigenvalues,
  does there exists a proper cone~$\bK$ for which~$A_i\, \bK \subset \bK$?   In general, the answer is negative already for~$m=1$:  a matrix with a Perron eigenvalue may not have an invariant cone~\cite{V}. 
  Now we make the following standard assumption: 
 \begin{assum}\label{a.10}
The semigroup~$\cS$ is {\em irreducible}, i.e.,  
its matrices do not share a common real nontrivial invariant subspace. 
\end{assum}
 We also use  the same concept for an arbitrary set of matrices~$\cA$. Clearly, an irreducible set 
generates an irreducible semigroup. In what follows, we always suppose that 
Assumption~\ref{a.10} is satisfied. Then 
the case of one matrix  becomes irrelevant since one matrix is always reducible for~$d\ge 3$.  
One then can conjecture that every {\em irreducible} Perron semigroup   
 leaves a cone invariant. This, however, is not true either:  
 
 \medskip 
 
\noindent \textbf{Fact~2.} {\em 
 The group~$SO(3)$ of  rotations of~$\re^3$ is Perron and does not have an invariant cone.}  

 \medskip 
 
Thus, there is a  gap between the Perron property and the existence of a common invariant cone. 
However,  this gap is not large.   Our main result,  Theorem~\ref{th.10} (Section~3), asserts  that an irreducible  Perron semigroup~$\cS$ 
 always has an invariant cone, provided that some generic  condition is satisfied. After we prove this theorem and derive some corollaries,  we attack the last problem: to 
 analyse the exceptional cases, when this condition is not met and a Perron semigroup does not 
have an invariant cone. 
 \medskip 
 
 \noindent \textbf{Problem.} {\em 
 Classify all irreducible Perron semigroups.}  
\smallskip 

\medskip 
 
 \noindent \textbf{Related works.} Spectral properties of multiplicative matrix groups and semigroups 
 have been studied in various contexts. The {\em unisingular groups}, i.e.,  matrix groups where each element 
 has~$1$ as an eigenvalue, are applied in algebraic geometry, Galois theory,  representation theory, etc. See~\cite{PS25, Z25}
 for recent results and an extensive bibliography. They are considered for matrices over various fields, 
including~$\re$ and~$\co$~\cite{BDKOR, BO}. Note that unisingular groups are not necessarily Perron since 
the eigenvalue~$1$ may not be leading. A closely related problem is a classification 
of semigroups of matrices with spectral radius~$1$. It draws much attention in the recent literature~\cite{OR, P, PV17}. For semigroups with an invariant cone, an overview is given in Section~2. 
\smallskip 

 \medskip

 The overall structure of this article is as follows. In Section~2 we give a short survey   of the 
 invariant cone problem.  Then, in  Section~3,  
we formulate the main results. Their proofs are given in Section~4. 
Then we will concentrate on classifying all Perron semigroups.
A simple case of nilpotent semigroups (which are also Perron by definition)
is considered in Section~5. 
Section~6 gives a complete classification of all Perron semigroups in law dimensions
($d=2,3,4$), while for higher dimensions only partial  results 
are obtained (Section~7). In particular, we construct several exceptional classes of Perron semigroups 
that do not have invariant cones. Section~8 contains a discussion and 
 several open problems.

  In what follows, we denote by~$\rho(A)$ the spectral radius of a matrix~$A$. 
  The vectors are denoted by bold letters, the scalars are by usual letters. 
 The group of isomenties (orthogonal transforms) of~$\re^n$ is denoted as usual by $O(n)$, 
  the group of proper isomenties (with determinant~$1$) is~$SO(n)$. Those groups 
  on a subspace~$\bL \subset \re^n$ are denoted by~$O(\bL)$ and~$SO(\bL)$ respectively. 
  We denote by~${\rm rank}_+ \cA$ the minimal positive rank of a family of 
  $d\times d$ matrices~$\cA$, i.e., the minimal rank of its nonzero elements. We usually assume that 
  a basis in~$\re^d$ is fixed and identify the matrices with the corresponding linear operators.

\bigskip 

\begin{center}
\large{\textbf{2. The invariant cone problem}}
\end{center}
\bigskip

The problem of deciding the  existence of a common invariant cone for a set of matrices was addressed in many works.  Semigroups with invariant cones 
have special asymptotic and spectral properties and some results on nonnegative matrices can be extended to them, see~\cite{BZ, J12, SW2} and references therein. Such semigroups have been applied in 
combinatorics and number theory~\cite{GP13, WZ}, in
neural networks~\cite{EMT},  chaos, dynamical system, and control~\cite{BV, FLS, GS}, 
in the study of random matrix products~\cite{H,  Per, P10, PJ13} and joint spectral characteristic of linear operators~\cite{GP13, Mej}.

Despite many applications, a little is known about the solution of 
the  invariant cone problem. For general finite sets of integer matrices,  
this problem  is algorithmically undecidable~\cite{P10c}. 
Nevertheless, it can admit efficient solutions for 
some classes of matrices or under extra assumptions. 
 The case of one matrix is done in~\cite{V}. 
Semigroups  of  $2\times 2$ matrices are studied in~\cite{EMT, RSS}. 
See also~\cite{FV} for families of commuting matrices,~\cite{RSS}  
for diagonal matrices, and~\cite{BZ, MP} for the algorithmic aspect of this problem. 
To the best of our knowledge, 
for  arbitrary finite sets of $d\times d$ matrices, 
the only known criterion of the existence of an invariant cone    is presented in~\cite{P10c}.
That criterion is rather theoretical since it reduces  the invariant 
cone problem to finding  the so-called $L_1$-spectral radius, 
which is also a hard task. Theorem~\ref{th.10} formulated in the next section 
gives another type 
of criterion.

\bigskip 

\begin{center}
\large{\textbf{3. The main results}}
\end{center}
\bigskip

The  condition for a matrix semigroup to possess an invariant cone 
will be formulated in terms of multiplicities of leading eigenvalues.  
\begin{defi}\label{d.10}
The   geometric index (or simply {\em index}) of a matrix is the sum of geometric multiplicities of its 
leading eigenvalues.  The index of a matrix semigroup is the minimal index of its elements. 
The algebraic index is defined similarly for algebraic multiplicities. 
\end{defi}
Equivalently, the  index of a matrix is the total number of Jordan blocks corresponding 
to all its leading eigenvalues. In geometrical terms, the index of~$A$ is equal to the maximal 
dimension of invariant subspace on which~$A$ is similar to an orthogonal operator multiplied by
the spectral radius~$\rho(A)$. The  word~{\em index} will be reserved for the geometric index; its  algebraic counterpart will be called explicitly~{\em algebraic index}, which is the total number of 
leading eigenvalues counting multiplicities.

\begin{remark}\label{r.18}
{\em For an irreducible semigroup~$\cS$, 
it is quite exceptional to have 
the index bigger than two. This requires, in particular,  that none of the matrices from~$\cS$  has a unique 
leading eigenvalue (counting multiplicities) or a unique pair of leading complex conjugate eigenvalues. 
If~$\cS$ is Perron semigroup, then the latter case is impossible. 
In applications, a standard assumption  for a  matrix semigroup is that 
at least one element is {\em asymptotically rank-one}, i.e., 
has a unique simple leading eigenvalue~\cite{GM, H}.   
So, a ``generic'' Perron semigroup has index one.  
Indeed, if a semigroup is generated by products of matrices from some set~$\cA$, then 
a small perturbation of one of the matrices of~$\cA$  makes its index (and hence, the index of the whole semigroup)
equal to one. 
An example of a bigger index is given by groups of orthogonal matrices.   However,  there are other 
constructions, see Section~7.  }
\end{remark}

\begin{theorem}\label{th.10}
An irreducible   Perron  semigroup possesses an invariant cone, 
unless its  index is bigger than or equal to~three. 
\end{theorem}
Thus, if matrices of an irreducible  Perron semigroup do not share a common invariant  cone, then 
all of them must have indices at least~$3$. This means that 
``generically'' each Perron  semigroup does possess an invariant cone. On the other hand,  each semigroup 
with an invariant cone is obviously Perron.  Therefore, 
Theorem~\ref{th.10} gives almost a criterion for the invariant cone problem. 

If every matrix from a semigroup~$\cS$ has an invarinat cone (the cones may be different for different matrices), then 
$\cS$ is Perron due to the Krein-Rutman theorem. Combining with Theorem~\ref{th.10}, we obtain

\begin{cor}\label{c.17}
If every matrix of an  irreducible semigroup~$\cS$ possesses an invariant cone,
then they all share a common invariant cone,  
unless the index of~$\cS$ is bigger than or equal to  three. 
\end{cor}

What about irreducible semigroups of index~$2$ ? 
The following result gives a comprehensive answer. It turns out 
that such semigroups are never Perron.  
\begin{theorem}\label{th.20}
There are no irreducible Perron  semigroups of index two. 
\end{theorem}
Since all semigroups leaving a cone invariant are Perron, 
we obtain 
\begin{cor}\label{c.30}
If an irreducible  matrix semigroup has an invariant cone, then its index 
cannot be equal to two. 
\end{cor}

The 3D rotation group~$SO(3)$ demonstrates that the index three  in Theorem~\ref{th.10} cannot be increased. 
This example   can be generalized to all odd dimensions~$d$ and to all subgroups of~$SO(d)$
by the following  fact
proved in the next section. 
\begin{prop}\label{p.10}
Let~$d\ge 3$ be an odd number; then an arbitrary irreducible subgroup of~$SO(d)$ (in particular,~$SO(d)$ itself) 
is Perron and does not have an invariant cone. 
\end{prop}
Not all irreducible Perron semigroups without an invariant cone consist of 
orthogonal matrices. In Section~7   we present  other constructions of index 3 and larger.  Moreover, for low dimensions it is possible to 
find all Perron semigropus. We collect those results in Theorem~\ref{th.30} below. 

Let us remember that any closed semigroup of orthogonal matrices is a group.
We say that a matrix simigroup~$\cS$ is  {\em proportional to a subgroup of~$SO(d)$} if  there is a 
basis in~$\re^d$, in which every element of~$\cS$ is a matrix from~$SO(d)$  
multiplied by a nonnegative scalar (that depends on the matrix). In other words,~$\cS \setminus\{0\}$
contains only non-singular matrices and becomes a subgroup of~$SO(d)$ after 
dividing all elements by their spectral radii and  taking a closure.  Similarly we define proportionality to a subgroup of~$O(d)$.

\begin{theorem}\label{th.30}
For~$d=2$ and $d=4$, all irreducible Perron semigroups are 
precisely those that have an invariant cone;  

For~$d=3$, every irreducible Perron semigroup~$\cS$ either 
has an invariant cone or is proportional to a subgroup of~$SO(3)$;

For every~$d\ge 5$ which is not  a power of two, 
there are irreducible Perron semigroups of algebraic index less than~$d$ that do not have an invariant cone. 
\end{theorem}
The case~$d=2$ is actually solved by~Theorem~\ref{th.10}, while the cases~$d=3$ and~$d=4$
are more complicated and require some preliminary work.  For~$d\ge 5$, we have no complete classification but 
only 
several constructions. Apart from the case $d=2^n$, there always exist 
irreducible Perron semigroups that have no invariant cone and not proportional to a subgroup of~$O(d)$. 
 The case~$d=2^n$ is unsolved and left as an open problem.

The proofs of Theorems~\ref{th.10} and~\ref{th.20} and of Proposition~\ref{p.10} are 
given in the next section.   Then in Section~5 we prove Theorem~\ref{th.30}, obtain more results on Perron semigroups and take   a precise look at exceptional cases for Theorem~\ref{th.10}.

\bigskip 

\newpage

\begin{center}
\large{\textbf{4. Proofs}}
\end{center}
\bigskip 

We first establish Theorem~\ref{th.10} by applying Theorem~\ref{th.20}. 
This means that in the proof of Theorem~\ref{th.10} we assume that 
the index of~$\cS$ is equal to one, because Theorem~\ref{th.20} excludes the 
case of  
 index two. 
Then we give a proof of Theorem~\ref{th.20}. 

We need several auxiliary notation. For a given matrix semigroup~$\cS$, 
we denote~$\cS_0 \, = \, \bigl\{\, tA\, : \ t\ > 0, \, A\in \cS\, \bigl\}$
and
$$
\cS_1 \, = \, \bigl\{\, [\rho(A)]^{-1}\, A\, : \ 
 A\in \cS, \, \rho(A) > 0\bigr\}\, \cup \,  \{ A\, : \ 
A\in \cS, \, \rho(A) = 0\}.
$$
 Thus, each matrix from the set~$\cS_1$ 
has the spectral radius equal to either zero or one.   Clearly, all matrices from~$\cS$
have Perron eigenvalues if and only if so do all  matrices from~$\cS_0$.  
  The same is for the existence of a common  invariant cone
and for the geometric and algebraic  indices. The same is true for~$\cS_1$. 
\bigskip 

\begin{center}
\textbf{4.1.  Proof of Theorem~\ref{th.10}}
\end{center}
\bigskip

\smallskip 

We begin with the following fact according to which 
the existence of a common invariant cone (different  from the entire space)  guarantees the existence of a proper common invariant  cone.  
\begin{lemma}\label{l.20}
If an irreducible  matrix semigroup has a convex  invariant 
cone different from the entire~$\re^d$, then this cone is proper. 
\end{lemma}
{\tt Proof.} Let~$\bK$ be a convex invariant cone of a semigroup~$\cS$. 
If ${\rm dim}\, \bK\, < \, d$, then the linear span of~$\bK$ is a common invariant subspace 
for~$\cS$, which contradicts the irreducibility. Thus,~$\bK$ is full-dimensional. 

Since the cone~$\bK$ does not coincide with~$\re^d$, its apex, which is at the origin, 
belongs to the boundary of~$\bK$. Hence, by the  convex separation theorem, 
the origin can be separated from~$\bK$ by means of 
a nonzero linear functional~$\ba$. Thus, 
$(\ba, \bx) \ge 0$ for all~$\bx\in \bK$. Hence, $\ba$ belongs to the 
dual cone~$\bK^*\, = \, \{\by \in \re^d: \ \inf_{\bx\in \bK} (\by, \bx) \, \ge \, 0\}$. 
Therefore,  
$K^* \ne \{0\}$. Note that~$\bK^*$ is an invariant cone for the transposed 
semigroup~$\cS^T\, = \, \{A^T: \ A\in \cS\}$. 
Since~$\cS^T$ is also irreducible, we conclude that~$\bK^*$
is full-dimensional. This implies that the primal cone~$\bK$ is pointed. 
Thus, $\bK$ is a proper cone. 

{\hfill $\Box$}
\medskip

\noindent {\tt Proof of Theorem~\ref{th.10}.}  Invoking Theorem~\ref{th.20} 
we can assume that the index of~$\cS$ is equal to one, i.e., the set  
$\cS_1$ contains a matrix~$A$ of index one. 
Consider two cases: $\rho(A) > 0$ (and hence, $\rho(A) = 1$) and~$\rho(A) = 0$. 

If~$\rho(A) =1$, then  the Perron eigenvalue of~$A$ is equal to~$1$ and  
 has a unique 
Jordan block. Denote this block by~$J$ and let~$r \ge 1$ be its size. The moduli of 
all other 
 eigenvalues of~$A$ are strictly smaller than one. Therefore, in a suitable basis~$\{\be_i\}_{i=1}^d$, we have    
\begin{equation}\label{eq.Ak}
A^k \quad = \quad 
\left(
\begin{array}{ccc}
J^k & {}& 0\\
{} & {}& {}\\
0 & {} & B_k
\end{array}
\right)\ , \quad k\in \n\, , 
\end{equation}
 where  $\, B_k \, \to \, 0$ as 
$k \to  \infty$. Here the block~$J^k$ has dimension~$r$ and 
the block~$B_k$ is of dimension~$d-r$. If~$r=1$, then $J^k=1$; if~$r\ge 2$, then 
\begin{equation}\label{eq.Jk}
J^k \quad = \quad 
\left(
\begin{array}{ccccc}
1 & k & \frac{k(k-1)}{2}& \cdots & k \choose r-1\\
0& 1 & k & \cdots & k \choose r-2\\
0 & 0 & 1 & \ddots & \cdots\\ 
\vdots  & \vdots & \vdots & \ddots & k\\ 
0 & 0 & 0 & \cdots & 1
\end{array} 
\right)
\end{equation}
Consider the set~$\cV \, = \, \{ X\, \be_1, \, X  \in \cS\, \}$. 
Its conical hull~$\bK\, = \, { cone}\, (\cV)$ is a common invariant cone for~$\cS$, possibly 
coinciding with~$\re^d$.  We wish to show  
that~$\bK$ is actually contained in some half-space, 
and hence, by Lemma~\ref{l.20}, it is proper. 
Let us prove that this half-space is~$\bH \, = \ \{\bx \in \re^d: \ x_r \ge 0\}$, where~$x_r$ denotes  
the $r$th coordinate of~$\bx$.  If this is not true, then there is~$X \in \cS$
such that~$(X\be_1)_r < 0$ which means~$X_{r1} < 0$.

Denote~$X_{r1} = - c, \, c> 0$, and estimate all the entries of the matrix~$A^kX$
for large~$k$. Firstly, all rows of~$A^k$ except for the first one 
have the norms~$O(k^{r-2})$. Hence, the same is true for the 
rows of the matrix~$A^kX$. Now take a look at the first row. 
For the first element, we have 
$(A^kX)_{11} \, = \, 
\sum_{j=1}^{r}  {k \choose j-1}  X_{j1} \, = \, -c\, {k \choose r-1} \, + \, \sum_{j=1}^{r-1} {k \choose j-1} X_{j1} \, = \, -c\, {k \choose r-1} \, + \, O(k^{r-2})  \quad {\rm as}\  k\to \infty.$ 
For all the other elements~$(A^kX)_{1s}, \, s=2, \ldots , d$, we obviously have 
$(A^kX)_{1s} \, = \,  O(k^{r-1})  \quad {\rm as}\  k\to \infty.$  Making a normalization 
we consider the matrix~$\, P_k \, = \, {k \choose r-1}^{-1} A^kX$. 
Thus,~$(P_k)_{11} \, = \, -\, c \, + \, O(1/k)$, all other elements of 
the first row are bounded uniformly in~$k$. The elements of all other rows 
are~$O(1/k)$. This yields that one eigenvalue of~$P_k$ tends to~$-c$, the others 
tend to zero as~$k\to \infty$. Thus, for large~$k$, the matrix~$P_k$ does not have a Perron eigenvalue, which is impossible since~$P_k\in \cS_0$. 
\smallskip 

It remains to consider the case~$\rho(A) = 0$.   
All eigenvalues of~$A$ are zeros, i.e., $A$ is a nilpotent matrix. 
The leading eigenvalue has a unique Jordan block, hence,  
in the Jordan basis, $A$ has~$d-1$ ones over the main diagonal and all 
other elements are zeros. Therefore, for the matrix $A^{d-1} = J^{d-1}$, we have 
$(A^{d-1})_{1d} = 1$ and all other entries are zeros. Consequently, 
for every matrix~$X\in \cS$, the product~$A^{d-1}X$ has the first  row equal to 
the $d$th row of~$X$ and the other rows  are zeros. This yields that  
the only nonzero eigenvalue of~$A^{d-1}X$ is equal to~$X_{d1}$. 
Now we define~$\cV \, = \, \, X\, \be_1, \, X  \in \cS$. 
If all the matrices~$A^{d-1}X, \, X\in \cS, \, k\in \n$, 
have Perron eigenvalues, then~$X_{d1} \ge 0$ for all~$X\in \cS$,
i.e., the $d$th coordinate of the vector~$X\be_1$ is nonnegative.  
 Therefore, 
$\cV$ lies at the half-space~$\{\bx \in \re^d: \ x_1 \ge 0\}$, which completes the proof. 

{\hfill $\Box$}
\medskip

\bigskip 

\begin{center}
\textbf{4.2.  Proofs of Theorem~\ref{th.20} and of Proposition~\ref{p.10}}
\end{center}
\bigskip

We begin with $2\times 2$ matrices. In this case, all Perron semigroups 
of  index two, including  reducible ones,  admit a simple description. 
\begin{lemma}\label{l.10}
If~$d=2$, then, for every Perron semigroup~$\cS$ 
of index two, we have either~$\cS_1\, = \, \{I\}$
or~$\cS_1\, = \, \{I, M\}$, where~$M$ is an arbitrary 
matrix similar to the mirror symmetry.  
\end{lemma}
{\tt Proof.} Every matrix from~$\cS$ has index at least~$2$ (by the assumption)
and at most~$2$ (since the dimension is~$2$), so, the index is equal to~$2$.  
None of the matrices from~$\cS$ has zero leading eigenvalue, otherwise its geometric multiplicity 
must be equal to two, and hence, the matrix is zero. 
Thus, each matrix from~$\cS_1$ has the Perron eigenvalue equal to one. 
If~$A \in \cS_1$ and ${\rm det}\, A > 0$, then the second eigenvalue is equal to one, hence,~$A=I$  due to 
the index condition. Thus, $I$ is the only element of~$\cS_1$
with positive determinant. If there are no other matrices, then~$\cS_1 = \{I\}$. 
  If there are at least two  other 
  matrices~$B,C\in \cS_1$, 
  then they both have negative determinants. 
Therefore,~${\rm det}\, BC > 0$, and hence, $BC\, =\, t\, I$ with some~$t > 0$. 
Thus,~$\cS_1$ contains  at most two 
matrices with negative determinants, and they  are proportional to 
a pair of mutually inverse matrices.  
Since,~$B$ has the second eigenvalue equal to~$-1$, it follows that 
 $B^2 = I$ and~$B$ is similar to a mirror symmetry.  
Thus,~$\cS\, = \, \{I, B\}$.

{\hfill $\Box$}
\medskip 

\begin{cor}\label{c.10}
For~$d=2$,  every Perron semigroup
with index two is reducible. 
\end{cor}

{\tt Proof of Theorem~\ref{th.20}.} Let~$\cS$ be a Perron matrix semigroup and 
a matrix~$A\in \cS_1$ has index two. It has a leading eigenvector with the 
eigenvalue~$1$ and another eigenvector with eigenvalue either~$1$ or~$-1$. 
All leading eigenvalues must be real, otherwise, the index is at least three. 
We assume that the second eigenvalue is also equal to one, otherwise we replace~$A$ by~$A^2$.
The two-dimensional plane spanned by those two eigenvectors 
is denoted by~$\bL$. The restriction of the operator~$A$ to~$\bL$ is 
the identity operator. 

Consider the case when the algebraic multiplicity of the  eigenvalue~$1$ 
is equal to two. In this case, all other eigenvalues are smaller than one by modulus. 
Hence, the limit~$\bar A\, = \, \lim_{k\to \infty} A^k$ is a projection of~$\re^d$
to~$\bL$. All matrices of the set~$\cR = \{\bar A B\, : \ B\in \cS\}$  have images on~$\bL$.
On the other hand, their indices are not smaller than two and they have Perron eigenvalues.  Hence, 
the restriction of~$\cR$ to~$\bL$ is a Perron semigroup 
of index two. Invoking now Corollary~\ref{c.10} we 
conclude that the matrices from~$\cR$ have a common one-dimensional 
invariant subspace~$\ell  \in \bL$.  
Denote by~$\bV$ the preimage of~$\ell$ by the projection~$\bar A$. 
If there is~$C\in \cS$ such that~$C\ell$ does not lie in~$\bV$, then~$\bar A C\ell$ 
 does not lie in~$\ell$. This is impossible since~$\bar A C \in \cR$. 
 We see that~$C\ell \subset \bV$ for all~$C\in \cS$. 
 Therefore, ${ span}\, \{C\ell : \ C\in \cS\}$ is a nontrivial common invariant subspace
 of~$\cS$. This contradicts  the irreducibility of~$\cS$, which completes the proof 
 in case when the eigenvalue~$1$ of the matrix~$A$ has multiplicity two. 
 
It remains to consider the case when  the algebraic multiplicity of the  eigenvalue~$1$ 
of the matrix~$A$ exceeds two. 
 This means that~$A$ has two Jordan blocks corresponding to the eigenvalue~$1$, at least one of them is nontrivial, and all other eigenvalues are smaller than one by modulus. 
 Let the sizes of blocks be~$n$ and~$m$, $n\ge m$. Thus, $n\ge 2$. 
 Take the Jordan basis~$\{\be_k\}_{k=1}^d$, where the first~$n$ vectors correspond to the first 
 block and the next~$m$ correspond to the second one.  Thus, $A\be_1 = \be_1$
 and~$A\be_{n+1} = \be_{n+1}$. 
The linear span of the set~$\{B\be_1: \ B\in \cS\}$ 
is an invariant subspace of~$\cS$, which is nontrivial,
since it contains~$\be_1 = A\be_1$. Hence, it coincides with~$\re^d$. 
On the other hand, the set of vectors with zero~$n$th coordinate is 
a hyperplane. Therefore, there is~$B\in \cS$ 
such that~$(B\be_1)_n \ne 0$ and so~$B_{1n} = c\ne 0$. 
Denote~$A^{(k)} = {k \choose n-1}^{-1}A^k$. 
If~$m<n$, then~$(A^{(k)})_{11} = c$ while the elements of 
all rows except for the first one tend to zero as~$k\to \infty$. 
Therefore, for large~$k$, the leading eigenvalue of~$A^{(k)}$ is close  to~$c$ and all
other eigenvalues are small. Hence, the matrix~$A^k \in \cS_0$  
has a unique leading eigenvalue and its index is equal to one, which contradicts  the 
assumption. If~$m=n$, all the rows of~$A^{(k)}$ except for
the first and, maybe,  the~$(n+1)$st ones tend to zero. 
For large~$k$, this implies that the matrix~$A^{(k)}$ has 
either one leading eigenvalue, or
two simple leading eigenvalues. In the former case, 
~$A^{(k)}$ has index one, and hence, so does~$A^k$, which contradicts the assumption. 
The latter case is considered above in the first part of the proof.

{\hfill $\Box$}
\medskip

{\tt Proof of Proposition~\ref{p.10}.} Consider an arbitrary matrix~$A\in SO(d)$. Since all its non-real eigenvalues 
are split into pairs of conjugates, their product is positive.
Among the real eigenvalues, the total number of negative ones is 
even, because~${\rm det}\, A \, = \, 1$. Since~$d$ is odd, it follows that 
there is a positive eigenvalue. It is leading, since 
all the eigenvalues are  equal by modulus. Thus, every irreducible subgroup of~$SO(d)$ is Perron. 
It remains so show that it does not have an invariant cone. We prove more: 
every set~$\cA$ of orthogonal matrices sharing a common invariant cone~$\bK$
share a common Perron eigenvector. Consider the compact set~$\bK_1$, which is 
the intersection of~$\bK$ with the  unit ball centered at the origin. 
Clearly,~$A \bK_1 \subset \bK_1$ for each matrix~$A\in \cA$ and actually
~$A \bK_1 = \bK_1$, because~$|{\rm det}\, A| = 1$ and hence,~$A$ preserves the volume. 
This implies that the center of gravity~$\bg$
of  the set~$\bK_1$ is a fixed point of~$A$. 
Thus, $\bg$ is a common eigenvector of all matrices from~$\cA$.

{\hfill $\Box$}
\medskip

\bigskip 

\begin{center}
\large{\textbf{5. Semigroups of nilpotent matrices}}
\end{center}
\bigskip 

We begin the anaysis of Perron semigroups with the  simple case
when the leading eigenvalues of all matrices of~$\cS$ are zeros, i.e.,  
all matrices of~$\cS$ are nilpotent. 
The $d$th power of  a nilpotent matrix is zero, which follows 
by considering  its Jordan form. Clearly, every  semigroup that  consists of nilpotent matrices is Perron. 
\begin{theorem}\label{th.45}
If all elements of a  martix semigroup~$\cS$ are nilpotent, then there is a basis in~$\re^d$ in which every 
matrix from~$\cS$ is upper triangular with zero diagonal. 
\end{theorem}
Note that the product of arbitrary~$d$ upper triangular matrices with zero diagonal
is equal to zero.  So, Theorem~\ref{th.45} implies that  if every element of a semigroup~$\cS$ is nilpotent, then  
the product of $d$ arbitrary elements (with repetitions permitted)
is zero. That is why we use the notation~{\em nilpotent semigroup}. 
The main corollary of Theorem~\ref{th.45} is that {\em a nilpotent semigroup 
always has a common kernel and hence, reducible}. We reformulate this fact as follows: 
\begin{cor}\label{c.40}
Every irreducible matrix semigroup contains an element with a positive spectral radius. 
\end{cor}
In the proof of Theorem~\ref{th.45} we use the notion of the {\em joint spectral radius} 
of a compact set of matrices~$\cA$:  
$$
\rho(\cA) \ = \ \lim_{n\to \infty} \max_{A_1, \ldots , A_n \in \cA}\|A_n\cdots A_1\|^{1/n}\, .
$$ 
This limit exists and does not depend on the matrix norm, see~\cite{J09, RS60} for this and other properties 
of the joint spectral radius.  The Berger-Wang theorem~\cite{BW92} asserts that 
$$
\rho(\cA) \ 
 = \ \limsup_{n\to \infty}
\max_{A_1, \ldots , A_n \in \cA}\rho \bigl(A_n\cdots A_1\bigr)^{1/n}.
$$ 

\smallskip

{\tt Proof of Theorem~\ref{th.45}.} It suffices to show that all matrices from~$\cS$ have a common kernel and then 
apply induction with respect to the dimension. Denote by~$\cS'$ the closure of the 
set of normalized matrices~$
\frac{1}{\|A\|}\, A\, , \   A\in \cS \setminus \{0\}$.  
This set is compact but not necessarily a semigroup. 
If the matrices from this set do not have a common kernel, then 
the value~$m\, = \, \min\limits_{\|x\|=1}\max\limits_{A \in \cS'}  \|A\bx\|$
is strictly positive. For every~$\bx\in \re^d$, there is a matrix~$A 
\in \cS'$ such that~$\|A\bx\| \ge m \|\bx\|$. 
Therefore, for every~$n$, there exists a sequence~$A_1, \ldots , A_n \in \cS'$
such that~$\|A_n\cdots A_1\|\, \ge \, m^n$. 
Hence,  $\rho(\cS') \ge m$. On the other hand, 
every product of matrices from~$\cS'$ has spectral radius zero
and by the Berger-Wang theorem~$\rho(\cS') = 0$.

{\hfill $\Box$}
\medskip 

Thus, all nilpotent semigroups are reducible and hence, they are out of our consideration. 
Nevertheless,  
Theorem~\ref{th.45} will be of use in later sections. We attack the problem of classifying irreducible 
Perron simigroups with the case of small dimension.

\bigskip 

\begin{center}
\large{\textbf{6. Perron semigroups of dimensions~$ d \le 4$}}
\end{center}
\bigskip 

We are going to prove the first two statements of Theorem~\ref{th.30}, the last one, 
which concerns arbitrary dimensions, will be established in the next section. 
We start with some preliminary work.

\begin{center}
\textbf{6.1. Auxiliary results}
\end{center}
\bigskip

It seems obvious that the algebraic index of a semigroup cannot exceed its minimal positive 
rank (denoted as~${\rm rank}_+ \cS$), which is the minimal rank of its nonzero elements. 
Indeed, for each matrix~$A\in \cS$, the total multiplicity of leading eigenvalues  cannot 
exceed~${\rm rank}\, A$, unless those eigenvalues are zeros. However, if~$\lambda_{\max} = 0$, 
then this is not true. Indeed, for nilpotent 
semigroups (Section~5), the algebraic index is equal to~$d$, while ${\rm rank}_+\cS \le d-1$.  
\begin{prop}\label{p.20}
The algebraic index of an irreducible semigroup~$\cS$ does not exceed~${\rm rank}_+\, \cS$. 
\end{prop}
{\tt Proof.} Denote~$q = {\rm rank}_+\, \cS$, choose a matrix~$R \in \cS$ of rank~$q$
and denote by~$\bL$ its image. Thus, the rank of each nonzero matrix from~$\cS$
is at least~$q$ and~${\rm rank}\, R\, = \, q$.  
Consider the  semigroup~$\cS_q = \{RA, \, A\in \cS\}$, whose matrices 
map $\re^d$ to~$\bL$.  If there is~$A \in \cS$ such that~$\lambda_{\max}(RA) \ne 0$, 
then the algebraic multiplicity  of this eigenvalue as at most~$q$, which completes the proof. 
Otherwise, if~$\lambda_{\max}(RA) = 0$ for all~$A \in \cS$, then all elements 
of the semigroup~$\cS_q|_{\bL}$ is nilpotent. Hence, by Theorem~\ref{th.45}, this semigroup 
has a common kernel~$\bU \subset \bL$. For arbitrary nonzero~$\bx \in \bU$, consider 
the space $\bV = {\rm span}\, \{A\bx, \ A \in \cS\}$. Clearly, $R\bV = \{0\}$ and  $\bV$ is a common invariant
subspace of~$\cS$. If~$\bV = \{0\}$, then~$\bx$ belongs to the common kernel of~$\cS$ which contradicts  
to the irreducibility. If~$\bV = \re^d$, then the image of~$R$ is zero (because~$R\bV = \{0\}$), which is impossible 
since ${\rm rank}\, R\, = \, q$. 

{\hfill $\Box$}
\medskip

\begin{prop}\label{p.25}
Let~$\cS$ be  an irreducible semigroup of $d\times d$ matrices and~$i$ be its algebraic index;  then 

\textbf{a)} Every matrix from~$\cS$ of rank less than~$i$ is nilpotent; 

\textbf{b)} If~$\, i \, > \, \frac{d}{2}$, then~$\cS$ does not contain nonzero nilpotent matrices. 
\end{prop}
{\tt Proof.} \textbf{a)} If~$\rho(A) > 0$, then the algebraic index of~$A$ does not exceed its rank, 
hence, $\, {\rm rank}\, A\, \ge \, i$, which is a contradiction. 

 \textbf{b)} Assume~$\cS$ contains a nonzero nilpotent matrix~$A$. 
 Let~$k$ be the biggest number such that~$A^k \ne 0$. Clearly, 
 ${\rm rank}\, A^k \, \le  \, \frac{d}{2}$, and hence~${\rm rank}\, A^k < i$.  
 Let~$\cS' = \{D \in \cS:  \ {\rm rank}\,  D \, < \, i\}$. This is a semigroup and~$A^k \in \cS'$. 
From item a) it follows that~$\cS'$ is nilpotent and by Theorem~\ref{th.45}, 
it has  a nontrivial common kernel~$\bV$. Denote by~$\bL$ the kernel of~$A^k$. 
Obviously, $\bV \subset \bL$. Take an arbitrary nonzero~$\bx \in \bV$. 
 The irreducibility of~$\cS$ implies that there exists~$Q \in \cS$ such that~$Q\bx \notin \bL$. 
Hence, $A^kQ\bx \ne 0$. On the other hand~${\rm rank}\, A^kQ \, \le \, {\rm rank}\, A^k \, < \, i$, therefore, 
 $A^kQ \in \cS'$  and so~$A^kQ\bx = 0$. The contradiction completes the proof.

{\hfill $\Box$}
\medskip

\begin{prop}\label{p.30}
An irreducible  semigroup~$\cS$ of~$d\times d$ matrices has an algebraic index~$d$
if and only if there is a basis in~$\re^d$ in which all nonzero matrices from~$\cS$
are proportional to orthogonal matrices. 
\end{prop}
{\tt Proof.} 
An arbitrary nonzero matrix~$A \in \cS$ has eigenvalues of equal modulus. 
If this modulus is zero, then~$A$ is nilpotent, which is impossible in view of part \textbf{b} of Proposition~\ref{p.25}. 
Hence, all nonzero elements of~$\cS$ are non-singular. 
Dividing each of them by the modulus of its determinant to  the power~$1/d$, we may assume 
that~$|{\rm det}\, A| = 1$ for all~$A \in \cS$. In this case,  
every~$A$ is diagonalizable with all eigenvalues equal to one by modulus. 
Now Theorem~2 from~\cite{PV17} implies that they are diagonalizable in a common basis, i.e., 
 there exists a basis in~$\re^d$ 
in which all matrices from~$\cS$ are orthogonal.

{\hfill $\Box$}
\medskip

Now we consider low-dimensional semigropus. Since the case~$d=2$ follows directly from Theorem~\ref{th.10}, 
we begin with~$d=3$.

\bigskip

\begin{center}
\textbf{6.2. The case} $\mathbf{d=3}$
\end{center}
\bigskip

{\tt Proof of Theorem~\ref{th.30} for $d=3$.} If~$\cS$ 
has no  invariant cones, then, by Theorem~\ref{th.10}, its index is at least~$3$. 
 Hence, all 
matrices have eigenvalues of equal modules.  In view of Proposition~\ref{p.30}, 
this means that in a suitable basis all nonzero elements of~$\cS$
become proportional to orthogonal matrices. Dividing each of them by its 
spectral radius we obtain a semigroup~$\cS_1 \subset O(3)$. 
If there is an element~$Q \in \cS_1$
with ${\rm det}\, Q = -1$, then its eigenvalues are~$\{1, 1, -1\}$. 
Indeed,~$Q$ must have a Perron eigenvalue~$1$, and the product of the two remaining eigenvalues 
is equal to~$-1$, hence, they cannot be complex conjugate. Thus,  $Q$ is an orthogonal  reflection across some 
hyperplane~$\bV \subset \re^3$. We see that 
each element of~$\cS_1$ is either a proper isometry (if the determinant is equal to~$1$) or a reflection.
Denote the set of those reflections by~$\cR$. 
For every proper isometry~$A \in \cA$, we have~$AQ \in \cR$  and 
therefore,~$(AQ)^2 = I$. Multiplying both parts by~$A^{-1}$, we 
obtain~$QAQ = A^{-1}$. If~$A$ is a rotation about a line~$\ell$, then 
its conjugate operator~$QAQ^{-1} = QAQ = A^{-1}$ is a rotation about~$Q\ell$. Hence, 
$Q\ell = \ell$. This is possible in two cases: either~$\ell$ is orthogonal 
to~$\bV$ or~$\ell \subset \bV$. 
In the former case~$A$ and $Q$ commute, hence~$A^2 = (AQ)^2 =I$, which means that~$A$
is a symmetry (rotation by~$180^{\circ}$) about~$\ell$. Denote this symmetry by~$F$. 
Thus, {\em every proper isometry from~$\cS$ is either~$F$ or 
a rotation with axis on~$\bV$}. If all those rotations have the same axis~$\ell$, then 
~$\ell$ is a common eigenvector for matrices from~$\cS$.
This is obvious for all rotations from~$\cS$, we need to show this 
for all reflections.  This is simple:~$\cR$ contains a reflection~$Q'$ 
with respect to a hyperplane~$\bV'\ne \bV$, then arguing as above we conclude that
$\ell \subset \bV'$ and hence,~$Q'\ell = \ell$. 
Thus,~$\ell$ is a common invariant subpace of~$\cS$, which is a contradiction. 
 
It remains to consider the last case:~$\cS$ contains at least two rotations~$A_1, A_2$ 
with different axes~$\ell_1, \ell_2 \subset \bV$. Since this must be true for 
all  reflections from~$\cR$, it follows that~$\cR = \{Q\}$
(since~$\bV$ is a unique plane containing~$\ell_1$ and~$\ell_2$). 
The operator~$A_1A_2A_1^{-1}$ is a rotation about the axis~$A_1\ell_2$. 
Therefore,~$A_1\ell_2 \subset \bV$ and so,~$A_1\bV \subset \bV$. 
This yields  that the plane~$\bV$ is respected by every proper isometry from~$\cS$. 
It is respected also  by~$Q$, hence,~$\bV$ is a common invariant subspace of~$\cS$, which is a 
contradiction.

{\hfill $\Box$}
\medskip

\bigskip

\begin{center}
\textbf{6.3. The case} $\mathbf{d=4}$
\end{center}
\bigskip 

We begin with the following technical result: 
\begin{lemma}\label{l.35}
Let~$\bL_1,\bL_2$ be two different linear hyperplanes in~$\re^d$, $\bE_1\subset \bL_1$ be 
a $(d-1)$-dimensional ellipsoid centered at the origin and an ellipsoid~$\bE_2 \subset \bL_2$
be its projection (possibly, non-orthogonal)  onto~$\bL_2$. 
Then there exists exactly one more projection that maps~$\bE_1$ to~$\bE_2$. Those two projections in a suitable 
basis have the form: 
 \begin{equation}\label{eq.p1p2}
 P_1 \bx \ = \ \bigl(0, x_2 - x_1, x_3, \ldots ,  x_d  \bigr)\ ; 
 \qquad P_2(\bx)\  = \  \bigl(0, x_2 + x_1, x_3, \ldots ,  x_d  \bigr)
 \end{equation} 
\end{lemma}
{\tt Proof.} Take the vectors~$\{\ba_i\}_{i=3}^{d}$ of the half-axes of the
ellipsoid~$\bE = \bE_1 \cap \bL_2$. Since~$\bE_2$ is a projection of~$\bE_1$ onto~$\bL_2$, we 
have~$\bE \, = \, \bE_1\cap \bE_2$. Let~$\ba_{1}$ be a point of the surface of~$\bE_1$ where 
the tangent hyperplane is parallel to~$\bL_2$ and similarly with~$\ba_2$. In the basis~$\{\ba_i\}_{i=1}^d$ 
the ellipsoids have the equations: 
$$
\bE_1\  = \ \Bigr\{\bx \in \re^d: 
\sum_{k\ne 2} x_k^2 = 1, \, x_2 = 0\Bigr\}\ ; \qquad 
 \bE_2\ = \ \Bigl\{\bx \in \re^d: 
\sum_{k\ne 1} x_k^2 = 1, \, x_1 = 0 \Bigr\}.
$$
A projection~$P$ onto~$\bL_2$ parallel to a unit vector~$\bb = (b_1, \ldots , b_d)$ is 
 given by the formula~$P\bx \, = \, \bigl(0\, , \, x_2 - \frac{b_2}{b_1}x_1\, , \, 
 \ldots ,  x_d - \frac{b_d}{b_1}x_1  \bigr)$. The relation~$P\bE_1 = \bE_2$
 means  that the equations 
 $\sum_{k \ne  2}^d x_k^2 = 1$ and $\, \frac{b_2^2}{b_1^2}\, x_1^2 \, + \, 
 \sum_{k = 3}^d \bigl( x_k - \frac{b_k}{b_1}x_1 \bigr)^2 = 1$ are equivalent.
 This happens precisely when~$b_1 = \pm b_2, \, b_3 = \ldots = b_d =0$. 
 Thus there are exactly two projections that map~$\bE_1$ to~$\bE_2$ and they have the 
 form~(\ref{eq.p1p2}).

{\hfill $\Box$}
\medskip

{\tt Proof of Theorem~\ref{th.30} for $d=4$.} 
If an irreducible  Perron semigroup~$\cS$ 
does not possess an invariant cone, then its index is at least~$3$. 
Hence, it is either~$3$ or~$4$. If the algebraic index is equal to~$4$, then 
 Proposition~\ref{p.30}
yields that in a suitable basis in~$\re^4$, all elements of~$\cS$ are proportional 
to orthogonal matrices.  After normalization it can be assumed that~$\cS$ is a subgroup of~$O(4)$. 
It was proved in~\cite[Theorem~3.4]{BO} that every group of~$4\times 4$ matrices 
in which every element possesses an eigenvalue~$1$ is reducible. 
Hence,  any irreducible subgroup of~$O(4)$ contains a matrix without 
eigenvalue~$1$, i.e., without  Perron eigenvalues. Thus, there are no irreducible Perron subgroups of~$O(4)$. 
It remains to consider the case when both the algebraic and geometric indices of~$\cS$
are equal to~$3$. 

Replacing each matrix~$A\in \cS$ by the pencil~$\{tA: \ t \ge 0\}$, we 
obtain a  Perron semigroup  without  an invariant cone. 
After replacing~$\cS$ by the closure of this semigroup it may be assumed 
that~$\cS$ is closed and positively homogeneous: $\lambda \cS = \cS$ for all~$\lambda \ge 0$. 
Due to Proposition~\ref{p.25} (part \textbf{b}), $\cS$ does not contain nonzero nilpotent matrices, hence, 
the spectral radii of all nonzero matrices in~$\cS$ are positive.  Moreover, 
$\cS$ does not contain matrices of rank less than three, otherwise, the index of such a matrix 
will also be less than~$3$. 
Consider the set of normalized matrices~$\cS_1 \, = \, \{[\rho(A)]^{-1} A: \ A\in \cS, A\ne 0 \}$. 
There exists a matrix~$D \in \cS_1$ with exactly three 
semisimple leading eigenvalues. 
Denote by~$\bL$ the three-dimensional hyperplane spanned by
the corresponding leading eigenvectors.  There is a subsequence  of powers of~$D$ that converges to 
an operator~$P$ of projection onto~$\bL$ parallel to some line~$\ell$. 
After a suitable change of coordinates we assume that~$\ell \perp \bL$. 
Thus,~$P \in \cS_1$ is an orthogonal projection onto~$\bL$. Consider the semigroup
$$
\cS_{P}\ = \ \Bigl\{PA\, : \quad A \in \cS  \Bigr\}. 
$$ 
Each element of this semigroup maps~$\re^4$ to~$\bL$. 
Since~$P^2 = P$, it follows that~$P \in \cS_{P}$. 
For every nonzero~$A \in \cA$, 
the operators~$PA$ is nonsingular on~$\bL$ since~${\rm rank}\, A\, \ge \, 3$. 
Thus, 
the restriction of~$\cS_P$ to~$\bL$ consists of nonsingular operators on~$\bL$ (apart from zero operator). 
Moreover, this restriction is irreducible on~$\bL$. Indeed, if there is 
a proper subspace~$\bL' \subset \bL$ such that~$PA\,\bL' \subset \bL'$ for all 
$A \in \cS$, then the subspace~$\tilde \bL = {\rm span}\, \{\ell, \bL'\} \, = \, 
P^{-1}\bL'$ is invariant for~$\cS$. To  show this consider arbitrary~$\bx \in \bL'$
and its orbit~$\{A\bx : \, A\in \cS\}$. If it is contained in~$\tilde \bL$, then 
its linear span is a prober  invariant subspace for~$\cS$. 
Otherwise, for some~$A' \in \cS$, we have~$A'\bx \notin \tilde \bL$. However, 
in this case~$PA'\bx \notin \bL'$, which is a contradiction. 
Thus, the restriction of~$\cS_{P}$ to~$\bL$ is an irreducible semigroup. 
It is contained in~$\cS$ and therefore, the index of each of its matrices is 
at least~$3$. Hence, it is equal to~$3$ because those operators act in the 
$3$-dimensional space~$\bL$. Thus, the restriction of~$\cS_P$ to~$\bL$ is isomorphic to a semigroup in~$\re^3$ of index~$3$, hence,
by Proposition~\ref{p.30}, in a suitable basis in~$\bL$, all its operators are proportional to proper isometries.
Dividing each of those operators by its determinant and removing the zero operator 
we obtain a subgroup~$\cS_{\bL}$ of~$SO(3)$.  

Denote by~$\bR$ the unit Euclinean sphere in~$\bL$.  
We have~$PA\bR =  \bR$ for all~$A \in \cS$. By the irreducibility, 
there are operators~$A \in \cS$ that do not leave~$\bL$ invariant. 
Denote by~$\cS'$ the set of such operators (this may not be a semigroup). 
For every~$A \in \cS'$, the set~$\bE = A\bR$ is a $3$-dimensional ellipsoid  which does not lie in~$\bL$
and~$P\bE$ is a sphere in~$\bL$. Multiplying~$A$ by a suitable positive constant, it can be assumed 
that~$P \bE = \bR$. Denote~$A\bL = \bM$ and observe that
the operator~$AP$, whose image is~$\bM$, has three leading semisimple eigenvalues 
with eigenvectors in~$\bM$. Therefore,~$AP\bM = \bM$ and there are powers~$(AP)^k$ 
that converge to a projector to~$\bM$ parallel to some line~$m$. 
Denote this projector by~$Q$. Thus, each element~$A \in \cS'$ generates at least one 
projector~$Q\in \cS$ onto the hyperplane~$A\,\bL$ such that~$A\bR = Q\,\bR$. 
  If all of those projectors are along lines parallel to~$\ell$, then 
  $\ell$ is a common invariant subspace for~$\cS$, which is impossible. 
  Hence,  for some~$A \in \cS'$, the corresponding projector~$Q: \, \bL\to \bM$ is not parallel to~$\ell$. 
  Therefore, the operator~$P^{-1}: \, \bL \to \bM$ is a projector different from~$Q$, since it is parallel to~$\ell$. 
  We see that there are two different projections that map~$\bR$ to~$\bE = Q\,\bR$. 
  By Lemma~\ref{l.35}, in a suitable coordinates they are given 
  by formula~(\ref{eq.p1p2}) with~$d=4$: 
 $$
 P^{-1}\bx\ = \ \bigl(0, x_2 - x_1, x_3,  x_4  \bigr)\ ; \qquad 
Q\, \bx \ = \  \bigl(0, x_2 + x_1, x_3,  x_4  \bigr)\, . 
$$ 
  Then~$PQ$ maps~$\bL = \{\bx \in \re^4: \, x_2=0\}$  to itself by the formula
  $\bx \, \mapsto \, (-x_1, 0, x_3, x_4)$. The 
 matrix of this operator on~$\bL$ is~$PQ \, = \, {\rm diag}(-1, 1,1)$
 with the determinant equal to~$-1$.  
 On the other hand, this operator belongs to~$\cS_P$ and therefore, $PQ \in SO(\bL)$,   which is a    
 contradiction.

{\hfill $\Box$}
\medskip

\bigskip 

\begin{center}
\large{\textbf{7. High dimensions: $\mathbf{d \ge  5}$}}
\end{center}
\bigskip 

We give the proof of the last part of Theorem~\ref{th.30}
by constructing two families of Perron semigroups of dimensions~$ d \ge 5$
that do not have  invariant cones. 

Theorem~\ref{th.10} asserts that Perron semigroups possess invariant cones 
apart from some exceptions. As we saw in Section~5, for~$d=2$ and~$d=4$, there are no 
exceptions at all, while for~$d=3$ the only exception  is given by Proposition~\ref{p.10}:  
this is the class of semigroups proportional to irreducible subgroups of~$SO(3)$. 
For~$d\ge 5$, there are more examples. Firstly, we present a class of Perron semigroups 
of odd  dimensions where  each semigroup  consists of singular matrices. 
The second class is  for all dimensions~$d$ that have at least one odd proper divisor.  
In both cases, the algebraic index is less than~$d$. This covers all dimensions~$d$ except for integer powers of two.

\bigskip

\begin{center}
\textbf{7.1. Odd~$d \ge 5$. Perron semigroups of singular matrices}
\end{center}
\bigskip

In the space~$\re^d$ with odd~$d$, we consider the following two subspaces
$$
\bL_0 \ = \ \bigl\{\bx \in \re^d:  \, x_1 = x_2 = 0\bigr\}, \qquad 
 \bL_1\ = \ \{\bx \in \re^d:  \, x_3 = x_4 = 0\}
 $$ 
and the following four operators~$P_i, S_i, \, i=0,1,$ acting in~$\re^d$: 
\begin{equation}\label{eq.PS}
\begin{array}{lcl}
P_0 \,\bx \ & = & \bigl(\,  0\, , \,  0\,  , \,  x_3 + x_1\,  , \,  x_4 + x_2\,  , \,  x_5\, ,  \ldots , \, x_{d}\, \bigr)\, ; \\
Q_0 \,\bx \ & = & \bigl(\,  0\,  , \,  0\,  , \,  x_3 - x_1\,  , \,  x_4 - x_2\,  , \,  x_5\,  ,  \ldots , \, x_{d}\, \bigr)\, ; \\
P_1 \,\bx \ & = & \bigl(\,  x_1 + x_3\,  ,  \,  x_2 + x_4\,  , \,  0\,  , \,  0\,  , \,  x_5\,  ,  \ldots , \, x_{d}\, \bigr)\, ;\\
Q_1 \,\bx \ & = & \bigl(\,  x_1 - x_3\,  ,  \,  x_2 - x_4\,  , \,  0\,  , \,  0\,  , \,  x_5 ,  \ldots , \, x_{d}\, \bigr)\, .\\
\end{array}
\end{equation}
We make several observations: 
\smallskip 

\noindent \textbf{1.} The operators~$P_i$ and~$Q_i$ are projections (non-orthogonal)
onto the subspace~$\bL_i, \, i = 0,1$. 
\smallskip 

\noindent \textbf{2}. On the subspace~$\bL_0$, the operators~$P_0P_1$
and~$Q_0Q_1$ are identities, while the operators~$P_0Q_1$ and~$Q_0P_1$
are proper isometries of~$\bL_0$ with the diagonal matrices~${\rm diag}\, (-1,-1,1, \ldots ,1)$. 

\smallskip 

\noindent \textbf{3.} The same is true on~$\bL_1$: the operators~$P_1P_0$
and~$Q_1Q_0$ are identities, while~$P_1Q_0 \, = \, Q_1P_0 \, = \, {\rm diag}\, (-1,-1,1)$
are proper isometries of~$\bL_1$. 

\smallskip 

\noindent \textbf{4.}  Let~$U_1  \in SO(\bL_1)$ be an arbitrary proper isometry of~$\bL_1$; 
then on the space~$\bL_0$, all the operators~$P_0U_1P_0^{-1}, \, 
Q_0U_1P_0^{-1}, \, P_0U_1Q_0^{-1},$ and $P_0U_1Q_0^{-1}$ are proper isometries 
of~$\bL_0$.   The same relations are true after the interchange of~$0$ and~$1$. 
\medskip 

This implies that the set 
\begin{equation}\label{eq.Aps}
\cS_{pq}\ = \ \Bigl\{U_iP_i\,  ,\  U_iQ_i\ : \qquad  U_i \in SO(\bL_i), \quad i=1,2 \ \Bigr\}\,  
\end{equation}
is a Perron semigroup. Indeed, each of its operators maps the unit ball of the space~$\bL_0$
either to itself or to the unit ball of~$\bL_1$, in both cases this is a proper isometry. 
The same holds for the unit ball of~$\bL_1$. 
Since~$U_i$ has a Perron eigenvector, it follows that both~$U_iP_i$ and~$U_iS_i$ 
have  the same Perron eigenvector. Thus, the semigroup~(\ref{eq.Aps}) 
is Perron and consists of singular matrices. 
\smallskip 

\begin{prop}\label{p.40}
The Perron semigroup~$\cS_{pq}$ defined by~(\ref{eq.Aps}) is irreducible, has algebraic index~$d-2$ and
does not possess an invariant cone.   
\end{prop}
{\tt Proof.}  
If~$\cS_{pq}$ has a nontrivial invariant subspace~$\bV \subset \re^d$, 
then for an arbitrary nonzero vector~$\bx \in \bV$, 
at least one of the projections~$\, P_0\,\bx, \, P_1\, \bx\, $ is nonzero. 
Suppose~$P_0\bx \ne \nill$; then the vectors~$U_0P_0\,\bx$ for all possible $U_0 \in SO(\bL_0)$, span the 
whole subspace~$\bL_0$. Hence,~$\bL_0 \subset \bV$ and 
consequently,~$\bL_0 = \bV$, otherwise 
$\bL_0$ would coincide with the whole space~$\re^d$. 
Thus,~$\bL_1 = P_1\bL_0$ is also in~$\bL_0$, which is a contradiction.
 It remains to show that~$\cS_{ps}$ does not have an invariant cone.  
For every~$\bv\in \bL_0$, there exists an operator~$U \in SO(\bL_0)$ for which 
$\bv$ is a unique leading  eigenvector corresponding to the eigenvalue~$\lambda_{\max}=1$.
Hence, if~$\cS_{ps}$ had an invariant cone~$\bK$, then this cone would contain 
one of the vectors~$\pm \bv$. Assume this is~$\bv$. 
There exists~$U_0 \in SO(\bL_0)$ such  that~$U_0\bv = -\bv$, hence, $ - \bv \in \bK$, 
which is impossible since~$\bK$ is a pointed cone.

{\hfill $\Box$}
\medskip

\bigskip

\begin{center}
\textbf{7.2. Perron semigroups with the block structure}
\end{center}
\bigskip

For given~$m\ge 2$ and odd~$p\ge 3$, we consider the following semigroup~$\cB_{p,m}$ of
$pm\times pm$-matrices. 
Every matrix~$A \in \cB_{p,m}$ is defined by 
arbitrary matrices~$U_1, \ldots , U_m$ from~$SO(p)$
and by a 
map~$g$ from the set~$\{1, \ldots , m\}$ to itself. 
The matrix~$A$ is partitioned to~$m^2$ blocks of size 
$p\times p$, the block~$(g(i), i)$ contains~$U_i, \, i=1, \ldots , m$, 
 all other blocks are zeros.  
\begin{theorem}\label{th.40}
For every~$m \ge 2$ and odd~$p\ge 3$, 
 the semigroup~$\cB_{p,m}$  is irreducible, Perron, has the algebraic   index~$p$, 
 and does not possess an invariant cone. 
\end{theorem}
{\tt Proof.} 
Consider an arbitrary matrix~$A \in \cB_{p,m}$ constructed by a map~$g$ and by 
some~$U_1, \ldots, U_m \in SO(p)$. 
Let~$G \ = \, \{a_i\}_{i=1}^m$  be the directed graph 
such that every vertex~$a_i$ has precisely one outgoing edge~$a_i \, \to \, a_{g(i)}$.
It is easily shown that all nonzero eigenvalues of~$A$
are eigenvalues of the products~$U_{i_n}\cdots U_{i_1}$ along the cycles 
$i_1\to \cdots \to i_n, \ i_{n+1} = i_1$, 
of the graph~$G$. All those products belong to~$SO(p)$ and hence, 
possess Perron eigenvalues, which are also eigenvalues of~$A$. 
Thus, all matrices from~$\cB_{p,m}$ have Perron eigenvalues. For the 
map~$g(i) = 1, \, i=1, \ldots , m$, all the 
matrices~$U_i$ are situated in the first block row of~$A$. 
Hence, $p$ eigenvalues of~$A$ coincide with the eigenvalues of~$U_1$, and all other eigenvalues 
are zeros. Therefore, the index of~$A$ is equal to~$p$. Now let us prove that~$\cB_{p,m}$
does not have an invariant cone. Take the cyclic permutation~$g$ such that~$g(i)=i+1, \, i=1, \ldots , m-1$ and~$g(m)=1$. For every vector~$\bv \in \re^{pm}$ we denote by 
$\bv_i\in \re^p$ its $i$th block of size~$p$. Thus,  $\bv \, = \, (\bv_1, \ldots , \bv_m)$. 
For every~$i$, there exists a matrix~$U_i \in SO(p)$ such that~$U_i\bv_i\, = \, - \bv_{i+1}$, where we set~$\bv_{m+1} = \bv_1$. Then those matrices and the map~$g$ define the matrix~$A \in \cS_{p,m}$, for which~$A\bv = - \bv$. We see that every vector~$\bv \in \re^{pm}$
can be mapped to the opposite vector~$-\bv$ by means of a suitable matrix from~$\cB_{p,m}$. 
Hence, $\cB_{p,m}$ has no invariant cones. 

To show that it is irreducible,  
it suffices to prove that for every nonzero vectors~$\bu, \bv \in \re^{pm}$, there is~$A\in \cB_{p,m}$
such that~$(\bu, A\bv) \ne 0$. Consider the same decomposition: 
 $\bv \, = \, (\bv_1, \ldots , \bv_m)$ and, respectively, $\bu
  \, = \, (\bu_1, \ldots , \bu_m)$. There are indices~$i, j$
for which~$\bv_i\ne 0$ and~$\bu_j\ne 0$. Consider an arbitrary permutation~$g$
such that~$g(i) = j$. Then choose~$U_k, \, k = 1, \ldots , m$
so that~$(\bu_j, U_i\bv_i) > 0$ and~$(\bu_{g(k)}, U_k\bv_k) \ge 0$ for all other~$k$. 
Then~$(\bu, A\bv) > 0$

{\hfill $\Box$}
\medskip

Note that in the proof of irreducibility in Theorem~\ref{th.40}, 
we actually use only those matrices from~$\cB_{p,m}$ generated 
by permutations~$g$. In the proof that~$\cB_{p,m}$ has no invariant cone, we use 
only cyclic permutations~$g$. Let~$\cQ_{p,m}$ be the set of matrices 
from~$\cB_{p,m}$ that are generated by permutations~$g$. Obviously, this is a semigroup. 
It is Perron, as a subset of~$\cB_{p,m}$.  Applying the above argument  demonstrates that
~$\cQ_{p,m}$ is irreducible  and does not 
possess an invariant cone. Finally note that all matrices from~$\cQ_{p,m}$ are orthogonal, 
This proves the following 
\begin{prop}\label{p.50}
For every~$m \ge 2$ and odd~$p\ge 3$, the set~$\cQ_{p,m}$ is a subgroup of~$O(pm)$. 
This subgroup is Perron, irreducible, and does not have an invariant cone. 
\end{prop}
This means that for every composite number~$d$ which is not a power of~$2$, 
there are irreducible Perron subgroups of the orthogonal group~$O(d)$ that contain at 
least one improper isometry. 

\bigskip

\begin{center}
\large{\textbf{8. Discussion and open problems}}
\end{center}
\bigskip 

A generic matrix semigroup has at least one matrix with 
a unique singular leading eigenvalue or a unique pair of complex leading eigenvalues.
In both cases the index is less than three and hence, the assumption 
of Theorem~\ref{th.10} is satisfied. It is an exciting challenge  to classify all exceptions. 
More generally: to classify all irreducible  
semigroups of index~$3$ and higher. We see that in small dimensions 
this problem is solvable. In general, however, it looks  complicated 
due to the existence  of several principally different families of counterexamples
presented above. One can construct superpositions of those families using block matrices.  
Therefore, the variety of exceptional families can be quite reach.  Nevertheless, 
some reasonable tasks can be formulated. First of all, the case of dimensions of~$d=2^n$
stays unclear.  
\medskip 

\noindent \textbf{Problem 1}. {\em Do irreducible Perron semigroups of dimension~$d=2^n$
always have invariant cones?}
\bigskip 

We know that the answer is affirmative for~$n=1, 2$. For subgroups of the orthogonal 
group, Proposition~\ref{p.40} concerns only composite dimensions.  In fact, it suffices to solve 
Problem~1 for~$n=3$. If there are examples of irreducible Perron semigroups in dimension~$8$ without invariant cones, 
then they can be extended to all dimension~$2^n, \, n\ge 4,$ by the block structure,
similarly to subsection~7.2. For~$d=8$, one could apply the idea from~\cite{BO}, where 
a group of matrices 
with eigenvalue~$1$ was constructed as follows. We consider the linear space~$\cV$ of~$3\times 3$ matrices with zero trace. 
 It has dimension~8. Every non-singular $3\times 3$ matrix~$A$ acts on~$\cV$ as a linear operator 
 $\ X \, \to \, A^{-1}XA, \ X\in \cV$. Those operators form an irreducible 
group in~$\re^8$, and each has an eigenvalue one~\cite{BO}. 
However, we can note that this group is not Perron. Indeed, 
every diagonal matrix~$A \, = \, {\rm diag}\, \{\lambda_1,\lambda_2,\lambda_3\}$, where
$|\lambda_1| > |\lambda_2| > |\lambda_3|$,  generates an operator with 
eigenvalues~$\bigl\{\frac{\lambda_i}{\lambda_j}: \ i, j=1,2,3 \bigr\}$
(the eigenvalue~$1$ has multiplicity two). If~$\lambda_1\lambda_3 < 0$, then 
the leading eigenvalue~$\frac{\lambda_1}{\lambda_3}$ is negative, hence, the group is not Perron. 

\bigskip 

\noindent \textbf{Problem 2}. {\em For which prime~$d$, there exist  Perron subgoups of~$O(d)$
which contain at least one matrix with determinant~$-1$ ? 
}

\bigskip

Finally, what can be said about exceptions of Theorem~\ref{th.10} of small indices? 
\smallskip 

\noindent \textbf{Problem 3}. {\em Characterize all irreducible Perron semigroups of index~$3$ that 
do not have invariant cones. }

 \end{document}